\documentclass[12pt]{article}
\usepackage{amsmath, amsfonts, amssymb, amsthm} 
\usepackage{makeidx}
\makeindex
\usepackage[colorlinks = true, linkcolor = blue]{hyperref}
\usepackage[nottoc]{tocbibind}
\newtheorem{theorem}{Theorem}[section]

\newtheorem{corollary}[theorem]{Corollary}

\def\R{\mathbb R} \def\Z{\mathbb Z} \def\C{\mathbb C} 

\def\N{\mathbb N}
\def\C{{\mathbb C}} 
\def\Q{\mathbb Q}
\def\P{\mathbb P}
\def\T{\mathbb T}

\def\supp{{{\rm supp}\,}}

\def\Ad{\hbox{\rm Ad}}

\def\Aut{\hbox {\rm Aut}}
\def\Inn{\hbox{\rm Inn}}
\def\Aff{\hbox {\rm Aff}}

\def\GL{\mathrm {GL}}
\def\SL{\mathrm {SL}}
\def\PSL{\mathrm {PSL}}
\def\Sp{\mathrm {Sp}}
\def\<{\,<\!}
\def\>{\!>\,}
\usepackage{graphics}
\usepackage{epsf}
\usepackage{epsfig}

\def\Aut{\hbox {\rm Aut}}

\title{Actions of automorphism groups of Lie groups}

\author{S.G. Dani}

\begin{document}

\maketitle

\section*{Introduction}

The aim of this article is to discuss the actions on a connected Lie group by subgroups of its 
automorphism group. The automorphism groups are themselves Lie groups (not necessarily 
connected) and the actions have, not surprisingly, played an important role in the study of various topics, including 
geometry, dynamics, ergodic theory, probability theory on Lie groups, etc.   

We begin in \S\,1 with results on the structure of the automorphism groups $\Aut (G)$, $G$ a Lie group, generalities about their subgroups, connected components, etc..
The automorphism group of a connected Lie group can be realised as a linear group via association with the corresponding automorphism of the Lie algebra $\frak G$ of $G$, and  \S\,2 is devoted to relating $\Aut (G)$ to algebraic subgroups, and 
more generally ``almost algebraic" subgroups,  of $\GL(\frak G)$; in particular the connected component of the identity in $\Aut (G)$ turns out to be almost algebraic and this has found considerable use in the study of various topics discussed in the subsequent sections. In \S\,3 we discuss orbits of various subgroups of $\Aut (G)$ on $G$, and especially conditions for them to be dense in $G$. In \S\,4 we discuss invariant 
measures and various 
ergodic-theoretic aspects of the actions of subgroups of $\Aut (G)$  on $G$, and in \S\,5 some of the topics in topological dynamics in which various special features of $\Aut (G)$ play a role. The action of $\Aut (G)$ on $G$ induces in a natural way an action on the space of probability measures on $G$ and properties of this action have played an important role in various questions in probability theory on Lie groups, in terms of extending certain aspects of classical probability theory to the Lie group setting. These are discussed in \S\,6. 

While our main focus will be on automorphism groups of connected Lie groups, at various points along the way generalisations to the case of more general locally compact groups are mentioned with references to the relevant literature.

\tableofcontents

\section{Structure of the automorphism groups} 

In this section we discuss various structural aspects of the automorphisms of connected Lie groups. 

\subsection{Preliminaries}\label{prelim}

Let $G$ be a connected Lie group and $\frak G$ be the Lie algebra of $G$. We denote 
by $\Aut (G)$ the group of all Lie automorphisms of $G$ equipped with its usual topology,
corresponding to uniform convergence on compact subsets (see  \cite{H} Ch.\,IX for some details). 
To each $\alpha \in \Aut (G)$ there  
corresponds a Lie automorphism $d\alpha$ of $\frak G$, the derivative of $\alpha$. 
We may view $\Aut (\frak G)$ in a natural way as a closed subgroup of 
$\GL(\frak G)$, the group of nonsingular linear transformations of $\frak G$, 
considered equipped with its usual topology. 
Consider the map $d:\Aut (G) \to \Aut (\frak {G})$, given by   $\alpha \mapsto d\alpha$, for 
all $\alpha \in \Aut (G)$. 
As $G$ is connected, $d\alpha$ uniquely determines $\alpha$, and thus the map  is injective.

When $G$ is simply connected  the map  is also surjective,
by Ado's theorem (cf. \cite{HN}, Theorem~7.4.1). When $G$ is not simply connected 
the map is in general not surjective; in this case $G$ has the form $\tilde G/\Lambda$, 
where $\tilde G$ is the universal covering group of $G$ and $\Lambda$ is a discrete 
subgroup contained in the center of $\tilde G$, and the image of $\Aut (G)$ in $\Aut (\frak G)$ 
consists of those elements for which  the corresponding automorphism of $\Aut (\tilde G)$ 
leaves the subgroup $\Lambda$ invariant (note that $\frak G$ may be viewed also, canonically, 
as the Lie algebra of $\tilde G$), and in particular it follows that it is a closed subgroup.
This implies in turn that $d$ as above is a topological isomorphism (see also  \cite{H}, Ch.\,IX, 
Theorem~1.2). 
We shall view $\Aut (G)$ as a subgroup of $\Aut (\frak G)$, and in turn $\GL(\frak G)$, via
the correspondence.  

As a closed subgroup of $\GL(G)$, $\Aut (G) $ is Lie group, and in particular a locally compact 
group; it may be  mentioned here that the question as to when the automorphism group of a general 
locally compact group is locally compact is analysed in \cite{Wu}. 

We note that $\Aut (\frak G)$ is a real algebraic subgroup of $\GL (\frak G)$, namely 
the group of $\R$-points of an algebraic subgroup of $\GL(\frak G \otimes \C)$, defined over $\R$. 
While $\Aut (G)$ is in general not a real algebraic subgroup (when viewed as a subgroup of $\GL (\frak G)$
as above), various subgroups of $\Aut (G)$ being algebraic subgroups plays an important 
role in various results discussed in the sequel.

\subsection{Subgroups of $\Aut (G)$}\label{subgroups}

We introduce here certain special classes of automorphisms which play an important role 
in the discussion in the following sections. 

Let $G$ be a connected Lie group. 
For any subgroup $H$ of $G$ let $\Aut_H(G)$ denote the subgroup of $\Aut (G)$ consisting of automorphisms 
leaving $H$ invariant; viz. $\{\tau \in \Aut (G)\mid \tau (H)=H\}$. A 
subgroup $H$ for which $\Aut_H(G)=\Aut (G)$ is  called a {\it characteristic} 
subgroup of $G$. Clearly, for any Lie group $G$ the center of $G$, the commutator subgroup $ {[G,G]}$, the (solvable) radical, the 
nilradical are some of the characteristic subgroups of $G$. Similarly, together with any characteristic 
subgroup, its closure, centraliser, normaliser 
etc.  are characteristic subgroups. 

For each $g\in G$ we get an automorphism $\sigma_g$
 of $G$ defined by $\sigma_g(x)=gxg^{-1}$ for all $x\in G$, called the {\it inner automorphism}  
corresponding to $g$. For a subgroup $H$ of $G$ we shall denote by $\Inn (H)$ the 
subgroup of $\Aut (G)$ defined by $\{\sigma_h\mid h\in H\}$. For any Lie subgroup $H$, 
$\Inn (H)$ is a Lie subgroup of $\Aut (G)$; it is a normal subgroup of $\Aut (G)$ when $H$
is a characteristic subgroup.  
In particular $\Inn (G)$ is a normal Lie subgroup of 
$\Aut (G)$; the group $\Aut (G)/\Inn (G)$ is known as the group of {\it outer automorphisms} of
$G$. We note that   $\Inn (G) $ may in general not be a closed subgroup of $\Aut (G)$; this is 
the case,
for example, for the semidirect 
product of $\R$ with $\C^2$ with respect to the action under which  $t\in \R$ acts by $(z_1,z_2)\mapsto 
(e^{i\alpha t}z_1, e^{i\beta t}z_2)$ for all $z_1, z_2\in \C$, with $\alpha, \beta $ fixed nonzero 
real numbers such that $\alpha/\beta $ is irrational.

Let $Z$ denote the center of $G$ and $\varphi: G\to Z$ be a (continuous) homomorphism of $G$
into $Z$. Let $\tau:G\to G$ be defined by $\tau (x)=x\varphi (x)$ for all $x\in G$. It is easy
to see that $\tau \in \Aut (G)$; we call it the {\it shear automorphism}\index{shear automorphism}, or more specifically
{\it isotropic shear automorphism} (as in \cite{D-inv}), associated with $\varphi$. We note
that for any continuous homomorphism $\varphi$ as above the subgroup $\overline{[G,G]}$ is contained in 
the kernel of $\varphi$, and hence the associated shear automorphism $\tau$ fixes $\overline{[G,G]}$ pointwise. In particular, if $\overline{[G,G]}=G$ then there are no nontrivial shear automorphisms. The shear automorphisms  form
a closed normal abelian subgroup of $\Aut (G)$, say $S$. 
We note that if $A=
\{\tau \in \Aut (G)\mid \tau (z)=z \hbox{ \rm for all } z\in Z\}$ then $\Aut (G)$ is the semidirect product of 
$A$ and $S$. 

 If $H$ is the subgroup of $G$ containing $\overline{[G,G]}$ and such that 
 $H/\overline{[G,G]}$ is the maximal compact subgroup of  $G/\overline{[G,G]}$, then $G/H$ is
 a vector group and 
the set of continuous homomorphisms $\varphi$ of $G$ into $Z$ such that $H\subset \ker \varphi$ has the 
natural structure of a vector space, and in turn the same holds for the corresponding set of shear automorphisms; the dimension of the vector space  
equals the product of the dimensions
of $G/H$ and $Z$.  When $G/\overline {[G,G]}$ 
is a vector group (viz. topologically isomorphic to $\R^n$ for some $n$) $S$ is a connected algebraic subgroup of $\Aut (G)$.  On the other hand, when $G/\overline{[G,G]}$ is compact then there are only countably many 
distinct continuous homomorphisms of $G$ into $Z$, and hence only countably many shear 
automorphisms. 

For a connected semisimple Lie group $G$, $\Inn (G) $ is a subgroup of finite index in $\Aut (G)$; in the case 
of a simply connected Lie group this follows from the corresponding statement for the associated Lie algebra
(cf. \cite{HN}, Theorem~5.5.14),  and the general case follows from the special case, since the inner automorphisms 
of the simply connected covering group factor to the original group $G$. On the other hand, for a connected nilpotent Lie group the group of outer automorphisms
is always of positive dimension (see \cite{Jac-outer}, Theorem~4). An example of a $3$-step simply connected nilpotent Lie 
group $G$ for which $\Aut (G)=\Inn (G) \cdot S$, where $S$ is the group of all shear automorphisms of $G$, 
is given in \cite{D-step}; in particular $\Aut (G)$ is nilpotent in this case; a larger class of connected 
nilpotent Lie groups 
for which $\Aut (G)$ is nilpotent  is also described in~\cite{D-step}. 

We recall here that the center of a connected Lie group $G$ is contained in a connected abelian Lie subgroup
of $G$ (see \cite{HN}, Theorem~14.2.8, or \cite{H}, Ch.\,16, Theorem~1.2). In particular the center is a compactly generated abelian group and hence  has a unique maximal compact subgroup; we shall denote it by $C$.
Being a compact abelian Lie subgroup, $C$ is in fact the cartesian product of a torus with a 
finite abelian group. As  the unique maximal compact subgroup of $G$, $C$ is $\Aut (G)$-invariant and hence we get a continuous homomorphism 
$q: \Aut (G) \to \Aut (C)$, by restriction of the automorphisms to $C$. Since the automorphism group of  a compact 
abelian group  is countable, it follows 
that $\Aut (G)/\ker q$, where $\ker q$ is the kernel of $q$, is a countable group. 

\subsection{Connected components of $\Aut (G)$}\label{conncomp}

As a closed subgroup of  $\GL(\frak G)$, $\Aut (G)$ is a Lie group with 
at most countably many connected components (see \cite{Hoc1} for an extensions of this
to not necessarily connected Lie groups). 
It is in general not connected; for $\R^n$, $n\geq 1$, the automorphism group, 
which is 
topologically isomorphic to $\GL (n,\R)$, has two connected components, while 
for the torus $\T^n\ (=\R^n/\Z^n)$, $n\geq 2$, the automorphism group is in fact  an 
  infinite (countable) discrete group. 

We shall 
denote by $\Aut^0 (G)$ the connected component of the identity in $\Aut (G)$; it is an 
open (and hence also closed) subgroup of $\Aut (G)$. We denote by $c(G)$ the group of connected 
components, namely $\Aut (G)/\Aut^0(G)$. 

For a connected semisimple Lie group, since $\Inn (G)$  is of finite index in $\Aut^0(G)$ (see 
\S\,\ref{subgroups}), $\Aut^0(G)= \Inn (G) $. 
Also, for these groups  the number  of  connected 
components is finite; the number is greater than one in many cases (see \cite{Jac} and \cite{Mur}
for details on  
the group of connected components of $\Aut (G)$ for simply connected semisimple groups 
$G$; see also  the recent paper \cite{Gun} where splitting of $\Aut (G)$ into the connected 
component and the component group is discussed).

We note that $\Aut^0(G)$ acts trivially on the 
unique maximal compact subgroup $C$ of the center of $G$, and hence when the image of 
the homomorphism $q:\Aut (G) \to \Aut (C)$, as in \S\,\ref{subgroups}, is infinite, $c(G)$ is infinite. 

When $C$ as above is the circle group, 
$\Aut (C)$ is of order  two, and the image of $q$ has at most two elements, but
nevertheless $c(G)$ can be infinite. A natural instance of this can be seen in the following:
Suppose $G$ has  a closed normal subgroup $H$ such that $G/H$ is a torus of positive 
dimension; then there is a unique minimal subgroup with the property and it is invariant under the action of $\Aut (G)$ and hence, 
modifying notation,  
we may assume $H$ to be $\Aut (G)$-invariant.  Then we have (countably) infinitely many 
(continuous) homomorphisms $\varphi: G/H\to C$ and for each such $\varphi$ we have an 
isotropic shear automorphism of $G$ (see \S~\ref{subgroups}).  It can be seen that the shear  automorphisms corresponding to distinct homomorphisms belong to distinct connected components of $\Aut (G)$, 
and hence $\Aut (G) $ has infinitely many connected components.  Thus $c(G)$ is infinite in this case also. 

The group $c(G)$ is finite if and only if $\Aut (G)$ is ``almost algebraic" as a subgroup of 
$\GL(\frak G)$; see \S \ref{aa} for a discussion on conditions for $\Aut (G)$ to be almost algebraic.  

\subsection{Locally isomorphic Lie groups}

Any connected Lie group is locally isomorphic to a unique (upto Lie isomorphism) simply connected Lie group, namely its universal covering group. Now let $G$ be  a simply connected Lie group and $Z$ be the 
center of $G$. Then all connected Lie groups 
locally isomorphic to $G$ are of the form $G/D$, where $D$ is a discrete subgroup of $Z$. Moreover, for two discrete subgroups $D_1$ and $D_2$ of $Z$ the Lie groups $G/D_1$ and $G/D_2$ are Lie isomorphic if and only if there exists a $\tau \in \Aut (G)$ such that $\tau (D_1)=D_2$. Thus the class of Lie groups
(viewed up to isomorphism of Lie groups) locally isomorphic to a given connected Lie group $G$ is in canonical one-one correspondence with the orbits of the action of $\Aut (G)$ on the class of discrete subgroups of it center $Z$, under 
the action induced by the $\Aut (G)$-action  on $Z$, by restriction of the automorphisms to $Z$. 

For $G=\R^n$, $n\geq 1$, all connected Lie groups isomorphic to $G$ are of the form $\R^m \times \T^{n-m}$, with $0\leq m \leq n$, thus $n+1$ of them altogether. For the  group $G$ of upper triangular $n\times n$ unipotent matrices ($n\geq 2$), which is a simply connected nilpotent Lie group, the center is one-dimensional and all its nontrivial discrete subgroups are infinite cyclic subgroups that are images of one another under automorphisms of $G$; thus in this case there are only two non-isomorphic Lie groups locally isomorphic to $G$ (including the simply connected one), independently of $n$. On the 
other hand there are simply connected nilpotent Lie groups $G$ for which $\Aut (G)$ is a unipotent group (when viewed as a subgroup of $\Aut (\frak G)$) (see \cite{D-step}); hence in this case the $\Aut (G)$-action on the center has uncountably 
many distinct orbits, and therefore there are uncountably many mutually non-isomorphic connected Lie groups that are locally isomorphic to $G$. 

Now let $G$ be a connected Lie group with discrete center, say $Z$. Then $Z$ is finitely generated (see \S\,\ref{subgroups}) and hence $Z$ has only countably many distinct subgroups. Considering the indices of the subgroups in $Z$ it can also be seen that when $Z$ is infinite there are infinitely many subgroups belonging to distinct orbits of the $\Aut (G)$-action on the class of subgroups; thus in this case there are countably infinitely many connected Lie groups locally isomorphic to $G$. This applies in particular  when $G$ is the universal covering group of $\SL(2,\R)$. When $Z$ is finite the 
number of Lie groups locally isomorphic to $G$ is finite, and at least equal to the number of prime divisors of the order of $Z$. For  simple simply connected Lie 
groups the orbits of the $\Aut (G)$-action on the class of subgroups of the center have been classified completely in \cite{GK}. 

\section{Almost algebraic automorphism groups}\label{algebraic}

Let $G$ be a connected Lie group and $\frak G$ be the Lie algebra of $G$. 
Let $\GL(\frak G)$ be realised as $\GL(n,\R)$, where $n$ is the dimension of $\frak G$, via 
a (vector space) basis of $\frak G$. For $g\in \GL(\frak G)$ let $g_{ij}$, $1\leq i,j \leq n$, denote the 
matrix entries of $g$ and let $\det g$ denote the determinant of $g$. 
A  subgroup $H$ of $\GL(\frak G)$ being a real algebraic group is equivalent to the condition 
that it  can be expressed
as the set of zeros (solutions) of a set of   polynomials in $g_{ij}$, $1\leq i,j \leq n$, and $(\det g)^{-1}$ as the variables;  in the present instance, 
the field being  $\R$, it suffices to consider a single polynomial in place of the set of polynomials.  

We call a subgroup of $\GL(\frak G)$ {\it almost algebraic}\index{almost algebraic subgroup} if it is of finite index in a real algebraic 
subgroup.  
A real algebraic subgroup is evidently a closed (Lie) subgroup of $\GL(\frak G)$, and since a 
connected Lie group admits no proper subgroups of finite index it follows that any almost algebraic
subgroup is an open subgroup of a real algebraic subgroup, and in particular it is a closed subgroup 
of $\GL (\frak G)$. 
It is known that any real algebraic subgroup has only finitely many connected components (see 
\cite{BT}, Corollary~14.5, for a precise result on the number of components),
and hence the same holds for any almost algebraic subgroup. 
\subsection{A decomposition}

Let $G$ be a connected Lie group and $\frak G$ the Lie algebra of $G$. 
In what  follows we view  $\Aut (G)$ as a subgroup of $\Aut (\frak G)$ and of $\GL(\frak G)$, via the
identification introduced earlier (\S\,\ref{prelim}).  We recall that any connected Lie group admits a maximal 
torus (a subgroup  topologically 
isomorphic to $\T^n=\R^n/\Z^n$ for some $n\geq 0$),  and any two maximal tori are conjugate 
to each other in $G$ (see \cite{HN}, Corollary~14.1.4).  

\begin{theorem}\label{decomp}{\rm (cf. \cite{D-Aut}, \cite{PW})}
Let $G$ be a connected Lie group with Lie algebra $\frak G$. Let $T$ be a maximal torus in $G$. 
Then there exists a closed connected
normal subgroup 
$H$ of $G$ such that the following conditions hold:

i) $\Inn (H)$  is an almost algebraic subgroup of $\Aut (\frak G)$, and 

ii) $\Aut (G) = \Inn (H) \Aut_T (G)$. 
\end{theorem}

The subgroup $H$ chosen in the proof is in fact invariant under all $\tau \in \Aut (G)$ and consequently 
$\Inn (H)$ is a normal subgroup of $\Aut (G)$. Since $\Inn (H)\cap \Aut_T(G)$ is trivial this further implies 
that $\Aut (G)$ is a semidirect product of the two subgroups. We may also mention here that the subgroup $H$ as chosen in \cite{D-Aut}  contains $[G,G]$, the commutator subgroup of $G$. 
When $\Inn (G)$ is an almost algebraic subgroup of $\Aut (\frak G)$, $H$ as in Theorem~\ref{decomp}
can be chosen to be $G$ itself; this applies in particular when $G$ is an almost algebraic 
subgroup of $\GL(n,\R)$ for some $n\geq 1$. 

\subsection{Almost algebraic subgroups of $\Aut (G)$}\label{aa}

In view of Theorem~\ref{decomp} $\Aut (G) $ is almost algebraic if and only if $\Aut_T(G)$ is 
almost algebraic for a (and hence any) maximal torus $T$ in $G$. In this respect 
we recall the following specific results. 

\begin{theorem}\label{thm:1.2}{\rm (\cite{D-Aut}, \cite{W}, \cite{PW})}
Let $G$ be a connected Lie group and $\frak G$ be the Lie algebra of $G$. Then the following statements hold. 

i) If the center of $G$ does not contain a compact subgroup of positive dimension, then 
$\Aut (G)$ is an almost algebraic subgroup of $\GL(\frak G)$. 

ii) $\Aut^0 (G)$ is an almost algebraic subgroup of $\GL(\frak G)$;

iii) if $R$ is the solvable radical of $G$ and $T$  is a maximal torus in $G$  then $\Aut (G)$ is almost algebraic if and 
only if the restrictions  of all automorphisms from $\Aut_T(G)$ to $T\cap R$ form a finite group of automorphisms of $T\cap R$; in particular, $\Aut (G)$ is almost algebraic if and only if $\Aut (R)$ is 
almost algebraic. 
\end{theorem}

Assertions (i) and (ii) above were deduced
from Theorem~\ref{decomp} in \cite{D-Aut}; (ii) was proved earlier by D. Wigner \cite{W},
and~(iii)  is due to  W.\,H. Previts 
and S.\,T. Wu \cite{PW}, where improved proofs were also given for (i) and (ii).

We note in particular that, in the light of Theorem~\ref{aa}(ii), $\Inn (G)$ is contained in an almost algebraic subgroup of $\Aut (G)$, namely $\Aut^0 (G)$. This turns out to be useful in various contexts on account of 
certain properties of actions of almost algebraic subgroups (see for example 
\S\S \ref{4.4} and \ref{5.1}). 

Let $G$ be a connected Lie group. Let $T$ be  the unique maximal torus contained in  the center of $G$.   
Then  $G/T$ is a connected Lie group whose center 
contains no nontrivial compact subgroup of positive dimension, and hence $\Aut (G/T) $ is an 
almost algebraic subgroup of $\Aut (\frak G')$, where $\frak G'$ is the Lie algebra of $G/T$. 
Each $\alpha \in \Aut (G)$ induces
an automorphism of $\Aut (G/T)$, say $\bar \alpha$.  If 
$\eta :\Aut (\frak G) \to \Aut (\frak G')$ is the canonical quotient homomorphism, defined by
$\eta ((\alpha)=\bar \alpha$ for all $\alpha \in \Aut (G)$,  then the image 
of any almost algebraic subgroup of $\Aut (\frak G)$ under $\eta$ is almost algebraic (as the map is a restriction
of a homomorphism of algebraic groups). In particular $\eta (\Aut^0(G))$ is an almost algebraic subgroup.

\subsection{Groups with $\Aut (G)$ almost algebraic}

For a class of connected Lie groups the following characterisation, incorporating a partial converse of Theorem~\ref{thm:1.2}(i) is proved  in \cite{CW}; 
there the issue is considered for all Lie groups 
with finitely many connected components, but we shall restrict here to when $G$ is connected;  
 (the general case involves some technicalities in its formulation). 

\begin{theorem}\label{chenwu}{\rm (\cite{CW})}
Let $G$ be a connected Lie group admitting a faithful finite-dimensional representation. Then
$\Aut (G)$ is almost algebraic if and only if the maximal torus contained in the 
center of $G$ is of dimension at most one, and it is also the maximal torus in the radical of $G$. 
\end{theorem}

Let  $G$ be as in Theorem~\ref{chenwu}  and $C$ be the maximal torus of the 
center of $G$ and suppose that $C$  is one-dimensional.  In the context of the examples of 
Lie groups with  $\Aut (G)$ having 
infinitely many distinct connected components discussed in \S \ref{conncomp} it may be noted that if $C$ 
 coincides with the maximal torus in the 
radical of $G$, then there does not exist a closed normal subgroup $H$ such that $G/H$ 
is a torus, as used in the argument there. 

It could happen that $\Aut (G)$ may be  almost algebraic (which in view of Theorem~\ref{aa}(ii) 
is equivalent to $c(G)$ being finite), for a connected Lie group $G$ (not admitting a faithful finite-dimensional representation) even if its center contains a torus of dimension 
exceeding $1$ (in fact of any given dimension). An example of this was indicated in \cite{D-Aut} (page
451) and has been discussed in detail in  \cite{PW}. The main idea involved in the example is that  Lie 
groups $G$ can be constructed such that  the central torus, though of higher dimension,
is the product of one-dimensional tori, each of which is invariant under a subgroup  of finite 
index in $\Aut (G)$. It is not clear whether there could be situations, with $c(G)$ finite, for which 
this may also fail, namely with no one-dimensional compact subgroups contained in 
the center and invariant under a subgroup of finite index in $\Aut (G)$. It is proved in \cite{PW} (Proposition~3.1), however
that for a connected Lie group $G$ of the form $H\times \T$, where $H$ is a connected Lie group and $\T$ is the one-dimensional 
torus, $\Aut (G)$ has infinitely many connected components, and hence is not almost algebraic.

\subsection{Automorphisms preserving additional structure}

When the Lie group has additional structure, the group of automorphisms preserving the structure would also be of interest. In this subsection we briefly recall some results in this respect. 

By a result of Hochschild and Mostow \cite{HM} if $G$ is a connected complex affine algebraic 
group then the connected component of the identity in the group $A$ of rational automorphisms of $G$ 
is algebraic (with respect to a canonical structure arising from the associated Hoff algebra), and 
moreover $A$ is itself  algebraic if either the center of $G$ is virtually unipotent (namely if it admits a subgroup of finite 
index consisting of unipotent elements) or the center of a (and hence any) maximal reductive subgroup of $G$ is of dimension at most $1$; the result may be compared with Theorem~\ref{thm:1.2} for real Lie groups. 
It may be noted that the results of \cite{HM} are in the framework of automorphism groups of 
affine algebraic groups over an algebraically closed field of characteristic zero. Some further elaboration 
on the theme  is provided by Dong Hoon Lee \cite{Lee}. 

For the connected component of the group of complex analytic automorphisms 
of a faithfully representable complex analytic group, a result analogous to 
that of  Hochschild and Mostow  
recalled above was proved by Chen and Wu~\cite{CW2}. 

\subsection{Linearisation of the $\Aut (G)$ action on $G$}\label{linrep}

Under certain conditions, a connected Lie group can be realised as a subset of $\R^n$, 
or the projective space $\P^{n-1} $, 
for some $n\geq 2$, in such a way that the automorphisms of $G$ are restrictions of linear
or, respectively, projective 
transformations. We call this {\it linearisation} of the $\Aut (G)$-action on $G$. We next 
discuss various results in this respect. 


Let $G$ be a connected Lie group. By an {\it affine automorphism}\index{affine automorphism} of $G$ we mean a transformation
of the form $T_g\circ \tau$ where $\tau \in \Aut (G)$,  and $T_g$ is a left translation by an element $g$
in $G$, namely $T_g(x)=gx$ for all $x\in G$. We denote by $\Aff (G)$ the group of all affine 
automorphisms. We identify $G$ canonically as a subgroup of $\Aff (G)$, identifying $g\in G$ 
with $T_g$ as above. Then $\Aff (G)$ is a semidirect product of $\Aut (G)$ with $G$, with 
$G$ as the normal subgroup, and we shall 
consider it  equipped with the Cartesian product topology.  We shall denote by $\Aff^0 (G)$ 
the connected component of the identity in $\Aff (G)$. 

Let $\frak A$ and $\frak B$ be the 
Lie algebras of $\Aut (G)$ and 
$\Aff (G)$ respectively. Let $\Ad : \Aff (G) \to \GL (\frak B)$ be the adjoint representation of $\Aff (G)$. 
Let $a$ be the dimensions of $\Aut (G)$ and  $V=\wedge^a\frak B$, 
the vector space of $a$-th 
exteriors over $\frak B$. Let  $\rho : \Aff (G) \to \GL(V)$ be the representation arising as the $a$-th exterior 
power of $\Ad$; we call 
$\rho$ the {\it linearising representation}\index{linearising representation} for $G$ (for reasons that would become clear below). 
Let  $L$ be the vector subspace $\wedge^a \frak A$ of $V=\wedge^a \frak B$; since $a$ is the dimension 
of $\frak A$, $L$ is a one-dimensional subspace. 

Consider the (linear) action of $\Aff (G)$ on $V$ via the representation $\rho$. Let $P=P(V)$ denote 
the corresponding projective space, consisting of all lines (one-dimensional subspaces) in $V$, equipped with its usual topology. The $\Aff (G)$-action on $V$ induces an action of $\Aff (G)$ on $P$. 
We have the corresponding 
actions of $\Aut (G)$ and $G$ by restriction of the action ($G$ being viewed as the group of 
translations as above). 

\subsection{Embedding of $G$ in a projective space}\label{linearisation}

Let $p_0$ denote  the point of $P$ corresponding to the  line $L$ as above. Let $S$ be the stabiliser of $p_0$ 
under the action of $G$. 
Consider any $g\in S$. Then $L$ is
invariant under the action of $g$ on $V$ and this means that the subspace $\frak A$ is invariant under the 
action of $g$ on $\frak B$. In turn we get that $\Aut^0 (G)$ is normalised by $g$ in $\Aff (G)$, 
and since $\Aff^0 (G)$ is a semidirect product of $\Aut^0 (G) $ and $G$ this implies that $g$ is 
contained in the center of $\Aff^0 (G)$.  Conversely it is easy to see that every element of the center 
of $\Aff^0 (G)$ fixes $p_0$. Thus $S$ coincides with  the center of $\Aff^0 (G)$.  We note also that $$S=\{g\in G\mid \tau (g)=g \hbox { \rm for all } 
\tau \in \Aut^0(G)\}.$$  

The orbit map $g\mapsto gp_0$ 
of $G$ into $P$ induces a canonical continuous bijection $j$ of $G/S$ onto its image in $P$, defined by $j(gS)=gp_0$ for all $g\in G$. When 
$S$ is trivial $j$ defines a continuous embedding of $G$ into $P$. We note that for any $\tau \in \Aut (G)$ 
and $g\in G$ we have $\tau (g)p_0=\tau g \tau^{-1}p_0=\tau g p_0=\tau (gp_0)$. 
Thus the map $j$ is equivariant with respect to the actions of $\Aut(G)$ on $G/S$ and $P$, viz. $j(\tau (gS))=\tau j(gS)$ for all $g\in G$ and $\tau \in \Aut (G)$. Thus when $S$ is trivial
the orbits of the $\Aut (G)$-action on $G/S$  are in one-one correspondence  in a natural way  with
orbits of the $\Aut (G)$-action on $P$ that are contained in the image of $G/S$.  

Recall that $\Aut^0(G)$ is an almost algebraic subgroup of $\GL(\frak G)$ (see Theorem~\ref{thm:1.2}(ii)). From this 
it can be deduced that the 
restriction of the representation $\rho$ as above to $\Aut^0(G)$ is an algebraic representation, 
viz. restriction of an algebraic homomorphism of algebraic groups. Hence every  orbit of 
$\Aut^0(G)$ on $P$ is locally closed, namely open in its closure (see \cite{Pr}, Lemma~1.22, for instance). 
Thus the above argument
shows the following. 

\begin{theorem}\label{locallyclosed}
Let $G$ be a connected Lie group and $S$ be the center of $\Aff^0(G)$. Then every orbit 
of the action of $\Aut^0(G)$ on $G/S$ is  locally closed (viz. open in its closure). 
\end{theorem}

The subgroup $S$, which  is noted to be  contained in the center of $G$, 
contains the maximal compact subgroup of the center. In particular 
if $S$ is trivial then  by Theorem~\ref{thm:1.2}(i) $\Aut^0(G)$ has finite index in $\Aut (G)$ and  hence the 
above theorem implies that all orbits of $\Aut (G)$ on $G$ are locally closed. 

We recall here that the condition of the orbits being locally closed is well studied in a wider context 
and is equivalent to a variety of other conditions of interest; see \cite{Gl}; 
(see also \cite{Ef} for a more general result). It represents in various 
ways the opposite extreme of the action being ergodic. 

\subsection{Embedding in a vector space}\label{linemb}

In analogy with the embedding of $G$ as a subset of $P$ (modulo the subgroup $S$ as above) 
we can also get an embedding of $G$ in the vector space $V$ as defined above, as follows. Let the notation be as
above and let $v_0$ be a nonzero point of $L$. We note that $S$ as above is also the stabilizer 
of $v_0$ under the $G$-action on $V$, since if $g$ fixes $v_0$ it fixes $p_0$, and if it fixes 
$p_0$ then we have $\tau (g)=\tau g \tau^{-1} =g$
for all  $\tau \in \Aut^0(G)$  and hence it fixes $v_0$. 
As $L$ is invariant under the action of $\Aut (G)$ 
we get a continuous homomorphism $s:\Aut (G) \to \R^*$ such that for all $\tau \in \Aut (G)$, $\tau (v_0) =s(\tau) v_0$. 
For any $g\in G$ and $\tau \in \Aut (G) $ we have $\tau (g)v_0=(\tau g \tau^{-1})v_0= \tau g (\tau^{-1}(v_0))=s(\tau^{-1} )\tau g v_0$. Consider the (linear) action of $\Aut^0(G)$ on $V$  such that $\tau \in \Aut^0(G)$ acts  by $v\in V \mapsto s(\tau^{-1})\tau v$. Then  the $\Aut^0(G)$-orbits on 
$G/S$ are in canonical correspondence with the orbits on $V$ under this  
$\Aut^0 (G)$-action. Under this action also the orbits are locally closed, by algebraic group considerations. 

\subsection{Algebraicity of stabilizers}

Let $G$ be a connected Lie group and $\rho :\Aff (G) \to \GL(V)$ be the corresponding linearising representation as in \S \ref{linrep}. We recall that $\Aut^0(G)$ is an almost algebraic subgroup of $\GL (\frak G)$, 
$\frak G$ being the Lie algebra of $G$, and that the restriction of $\rho$ to $\Aut^0 (G)$ is the restriction of a homomorphism of algebraic groups. In particular, for any $v\in V$ the stabiliser $\{\tau \in \Aut^0(G) \mid \tau (v)=v\}$ is an almost algebraic subgroup of $\Aut^0(G)$ (algebraic subgroup if $\Aut^0(G)$ is algebraic). Let $v_0$ be as in  \S \ref{linemb} and $S$, as before, be the center of $\Aff^0 (G)$, which is the stabiliser of $v_0$ under the action as in \S \ref{linemb}. Then by the preceding observation 
for all $g\in G$ the subgroup $\{ \tau\in \Aut^0(G) \mid \tau (gS)=gS\}$ is  almost algebraic. In particular for a connected Lie group $G$ such that the center of $\Aff^0(G)$ is trivial, for all $g\in G$ the stabiliser 
$\{ \tau\in \Aut^0(G) \mid \tau (g)=g\}$ of $g$ under the action of $\Aut^0(G)$ is an almost algebraic subgroup;
this holds in particular when the center of $G$ is trivial. 

Let $G$ be a connected Lie group and for $g\in G$ let $S(g)$ denote the stabiliser $\{\tau \in \Aut (G)\mid \tau (g) =g\} $.  It is shown in \cite{Dj} that when $G$  is a simply connected solvable Lie group, $S(g)$ is  
 an algebraic subgroup for all $g\in G$.  In the general case it is shown that the connected 
 component of the identity in $S(g)$ is an almost algebraic subgroup for all $g$ which are of the form $\exp \xi$ for 
 some $\xi $ in the Lie algebra of $G$, $\exp$ being the exponential map associated with $G$. It is noted that 
 $S(g)$  itself need not be an algebraic subgroup, as may be seen in the case when  $G$ is 
the universal covering group of $\SL (2,\R)$  and $g$ is one of the generators of the center of $G$ (the latter is an infinite cyclic subgroup); in this case $S(g)=\Inn (G)$, and it is of index $2$ in $\Aut (G)$, but not an algebraic subgroup.  

For elements $z$ contained in  the center of $G$ it is proved in \cite{Go2} (see also~\cite{D-Aut}) that 
 the connected component of the identity in $S(z)$ 
 is an almost algebraic subgroup, and if $G$ has no compact central subgroup of positive dimension then 
 $S(z)$ itself is almost algebraic. It is however not true that $S(z)$ is almost algebraic whenever $\Aut (G)$ 
is almost algebraic; a counterexample in this respect may be found in \cite{PW} (p. 432).  

We note in particular that if  $A$ is the subgroup of $\Aut (G)$  consisting of 
all automorphisms which fix the center pointwise then the connected component of the identity 
in $A$ is an almost algebraic subgroup, and if $G$ has no compact central subgroup of positive 
dimension then $A$ is almost algebraic.

\section{Dense orbits}

It is well-known that an automorphism $\alpha$ of the $n$-dimensional torus $\T^n$, $n\geq 2$,
admits dense orbits (and is ergodic with respect to the Haar measure; see \S \ref{erg} for a discussion 
on the ergodicity condition) when $d\alpha$
has no eigenvalue which is a root of unity.  
These  automorphisms constitute some of the basic examples in 
ergodic theory and topological dynamics, and have been much studied for 
detailed properties from the point of view of the topics mentioned (see \cite{CFS}, \cite{Wal},
for instance; see also \cite{Sch} and \cite{D-encl} for generalisations). 
In this section we discuss actions of subgroups of $\Aut (G)$ with a dense 
orbit on $G$. 

\subsection{$\Aut (G)$-actions with dense orbits}\label{denseorb}

Let $G$ be a connected abelian Lie group.  Then it has the form $\R^m\times \T^n$ for
some $m,n\geq 0$ and in this case it is easy to see that the $\Aut (G)$-action on 
$G$ has dense orbits, except when $(m, n)=(0,1)$ (viz. when $G$ is the circle group); if $n=0$ 
the complement of $0$ in $\R^m$ is a single orbit. 
If $m\geq 1$ then the action of $\Aut^0(G)$  also has a dense orbit on $G$, as
may be seen using the isotropic shear automorphisms (see \S \ref{subgroups}) associated with homomorphisms
of $\R^m$ into $\T^n$. 

Next let $G$ be a $2$-step connected nilpotent Lie group\index{two-step nilpotent group}, namely $[G,G]$ is contained in the 
center of $G$. Let $Z$ denote the center of $G$. Then $G/Z$ is simply connected,
and hence is Lie isomorphic to $\R^n$ for some $n$. The $\Aut (G)$-action on $G$
factors to an action on $G/Z$. Using shear automorphisms (see \S \ref{subgroups}) it can be seen 
that the $\Aut (G)$-action on $G$ 
has a dense orbit on $G$ if and only if  the $\Aut (G)$-action on $G/Z$ has a dense 
orbit. The latter condition  holds in particular if $G$ is a free 2-step simply connected Lie group; it may be
recalled that $G$ is called a free 2-step nilpotent Lie group if its Lie algebra is of the 
form $V\oplus \wedge^2V$, with $V$ a finite-dimensional vector space over $\R$, and 
the Lie product is generated by the relations $[u,v]=u\wedge v$ for all $u,v\in V$, and
$[u,v\wedge w]=0$ for all $u,v,w\in V$. When $G$ is a free 2-step simply connected 
nilpotent Lie group the quotient $G/Z$ as above corresponds to the vector space $V$ 
and every nonsingular automorphism of $V$ is the factor of a $\tau \in \Aut (G)$, which 
leads to the observation as above. 

For a general simply connected $2$-step nilpotent Lie group $G$ the Lie algebra $\frak G$ 
of $G$ can be expressed as $(V\oplus \wedge^2 V)/W$, where $V$ is the  vector space 
$G/Z$  ($Z$ being the center of $G$) and $W$ is a vector subspace of $\wedge^2 V$;  $V\oplus 
\wedge^2 V$ is the free $2$-step nilpotent  Lie algebra with the structure as above, and any vector subspace of 
$ \wedge^2 V$ is a Lie ideal in the Lie algebra (see \cite{AS}, where the example is discussed in a different  context). 
It is easy to see that in this case the image of $\Aut (G)$ in $\GL (V)$, under the map associating to 
each automorphism its factor on $V=G/Z$, is the subgroup, say $I(W)$,  consisting of $g\in \GL (V)$ such that the corresponding 
exterior transformation $\wedge^2 (g)$ of $\wedge^2 V$ leaves the subspace $W$ invariant.  
It can be seen that $I(W)$ has an open dense orbit on $V$ for various choices of $W$, and in these cases 
by the argument as above  $\Aut (G)$ has open dense orbit on $G$. For example, if $e_1, \dots, e_n$, $n\geq 2$ is a linear basis of $V$, then this is readily seen to hold if $W$ is the subspace spanned by $e_1\wedge e_2$ or, more generally,  by sets of the form $\{e_1\wedge e_2, e_3\wedge e_4, \dots, e_{2k-1}\wedge e_{2k}\}$, for $n\geq 2k$.  

Along the lines of the above arguments it can be seen that for any $k\geq 2$ if $G$ is a free $k$-step 
simply connected Lie group then the action of $\Aut (G)$ on $G$ has an open dense orbit. Also, given 
a simply connected $k$-step nilpotent Lie group $H$ there exists simply connected $k+1$-step nilpotent Lie group $G$ such that $H$ is Lie isomorphic to $G/Z$, where $Z$ is the center of $G$, and every automorphism of $H$  is a factor of an 
automorphism of $G$.  Therefore using the examples as above one can also construct examples, for 
any $k\geq 2$, of simply connected $k$-step nilpotent Lie groups such that the $\Aut (G)$-action on $G$ 
has an open dense orbit.

In the converse direction we have the following.  

\begin{theorem}\label{denseorbit} 
Let $G$ be a connected Lie group. Suppose that there exists $g\in G$ such 
that the closure of the $\Aut (G)$-orbit of $g$ has positive Haar measure in $G$. 
Then $G$ is a nilpotent Lie group. 
\end{theorem}

A weaker form of this, in which it was assumed that the orbit is dense in $G$, was proved in 
\cite{D-Sank} (Theorem~2.1), but the same argument is readily seen to yield the stronger assertion as above. 
By  a process of approximation, Theorem~2.1 in \cite{D-Sank} was extended 
 to all finite-dimensional connected locally compact groups, and in the same way one can also 
 get that the assertion in Theorem~\ref{denseorbit} holds also for all finite-dimensional
 connected locally compact groups. 

A nilpotent Lie group $G$ need not always have dense orbits under the 
action of $\Aut (G)$. Examples of $2$-step simply connected nilpotent Lie 
groups with no dense orbits under the $\Aut (G)$-action 
are exhibited in \cite{D-step}. There are also examples in \cite{D-step} of $3$-step simply connected 
nilpotent Lie groups for which $\Aut (G)$ is a unipotent group, namely when $\Aut (G) $ is viewed  a subgroup of $\GL (\frak G)$ all its elements are  unipotent linear transformation, and consequently all
orbits of $\Aut (G)$ on $G$ are closed, and hence lower dimensional, submanifolds of $G$. 

There has been a detailed study of groups in which $\Aut (G)$ has only finitely many, or countably
many, orbits, in the broader context of locally compact groups, and also abstract groups. The results in this respect for connected
Lie groups will be  discussed in~\S \ref{forbits}.

\subsection{Connected subgroups of $\Aut (G)$ with dense orbits}

Theorem~\ref{denseorbit} implies in particular that there is no connected Lie group $G$ for which 
$\Aut (G)$ has a one-parameter subgroup whose action on $G$ has a dense 
orbit; by the theorem, a Lie group $G$ with that property would be nilpotent and since it has to be noncompact,  going to a quotient we 
get that  for some $n\geq 1$, $\GL(n,\R)$ has a one-parameter subgroup
acting  on $\R^n$ with a dense orbit, but simple 
considerations from linear algebra rule this out. 

On $\R^{2n}$ viewed as $\C^n$, $n\geq 2$, we have 
linear actions of $ \C^{n-1}\times\C^{*}$,  which have an open dense orbit: the action of  $(z_1, \dots , z_{n-1}, z)$, 
where  $z_1, \dots,  z_{n-1} \in \C$  and $z\in \C^*$, is defined on the standard basis vectors $\{e_1, \dots, e_n\}$ by $e_j\mapsto ze_j$ for $j=1,\dots, n-1$ and $e_n\mapsto  \sum_{j=1}^{n-1}z_je_j + ze_n$.   The same also holds 
for (outer) cartesian products of such actions, and in particular we have linear actions of  $\C^{*n}$ on $\R^{2n}$,  with 
dense orbits. Since for $n\geq 2$ we can realise $\R^{n+1}$  as a dense subgroup of $\C^{*n}$, we get linear actions $\R^{n+1}$ on $\R^{2n}$ admitting dense (but not open) orbits; we note that in analogy with the above we can get linear actions of $\R^{*n}$ on $\R^n$ with open dense orbits, but they do not yield actions of  $\R^n$
with dense orbits. It can be seen using the Jordan canonical form that there is no linear action of $\R^2$ admitting dense orbits, and hence there is no action of $\R^2$ on any Lie group admitting  dense orbits, by an argument as above, using Theorem~\ref{denseorbit}.

On $\R^n$, apart from $\GL(n,\R)$, various proper subgroups act transitively on the complement
of $\{0\}$; e.g. $S\! \cdot \! O(n,\R)$ where $S$ is the subgroup consisting of nonzero scalar matrices
and $O(n,\R)$ is the orthogonal group, or the symplectic group $\Sp(n, \R)$ for even $n$. Also 
if $O(p,q)$ is the orthogonal group of a quadratic form of signature $(p,q)$ on $\R^n$, $n=p+q$
then the action of $S\!\cdot \!O(p,q)$, with $S$ as above, on $\R^n$ has an open dense orbit (the 
complement consists of  the set of zeros of the quadratic form, which is a proper algebraic subvariety of $\R^n$). 

Subgroups of $\GL(n,\R)$ acting with an open dense orbit have been a subject of much 
interest in another context; $\R^n$ together with such a subgroup is called a pre-homogeneous 
vector space\index{pre-homogeneous vector space}; the reader is referred to \cite{K} for details. 

\subsection{Automorphisms with dense orbits}

Analysis of the issue of dense orbits was inspired by a question raised
by P.\,R. Halmos in his classic book on Ergodic Theory (\cite{Hal}, page 29)\index{Halmos' question} as to whether a 
noncompact locally compact group can admit a (continuous) automorphism
which is ergodic with respect to the Haar measure of the group, namely such that there is no measurable set invariant under the automorphism such that both the set and its complement have positive Haar measure (see \S\,4.1 for more on ergodicity). An automorphism as in 
Halmos' question would have a dense orbit (assuming the group to be 
second countable), and in particular one may ask 
whether there exists an automorphism of a noncompact locally compact group with a 
dense orbit.  This question was answered in the negative in a paper of 
R. Kaufman and M.~Rajagopalan \cite{KR} and T.\,S. Wu \cite{Wu-auto} for connected locally compact groups, and 
N. Aoki \cite{Ao} for a general locally compact group (see also \cite{PW2} for some clarifications 
on the proof in \cite{Ao}); it may also be mentioned here that 
the analogous question is studied for affine automorphisms of locally compact 
groups (namely transformations which are  composite
of a continuous group automorphism with a translation by a group element) in  
\cite{D-dense} for connected groups, and in \cite{Kas} in the generality of  all locally compact groups. 

Though there are nonabelian compact groups admitting automorphisms with 
dense orbits, such as $C^{\Z}$ where $C$ is a compact nonabelian group, for which
the shift automorphism has dense orbits, a compact connected Lie group admits 
such an automorphism
only if it is a torus of dimension $n\geq 2$. For a compact semisimple Lie group $G$,
$\Inn (G) $ is a subgroup of finite index in $\Aut (G)$ and the orbits of $\Aut (G)$ on $G$ 
are closed submanifolds of dimension less than that of $G$, and in particular  not dense in $G$; 
a general compact
connected Lie group has a simple Lie group as a factor and hence the preceding 
conclusion holds in this generality also. 

\subsection{$\Z^d$-actions}

The analogue of Halmos' question for $\Z^d$-actions, namely the multi-parameter
case with $d\geq 1$, and more generally actions of abelian groups of automorphisms, was considered in \cite{D-Zd}, where the following is proved. 

\begin{theorem}\label{Zd}
Let $G$ be a connected Lie group. Suppose that there exists an abelian subgroup $H$ of $\Aut ( G)$  such that the $H$-action on $G$ has a dense orbit. Then there exists a compact subgroup $C$ contained in the center of $G$ such that $G/C$ is topologically isomorphic to $\R^n$ for some $n\geq 0$; in particular $G$ is a two-step nilpotent Lie group. If moreover the $H$-action leaves invariant the Haar measure on $G$ then $G$ is a torus. 
\end{theorem}

Recall that the $n$-dimensional torus $\T^n$, where $n\geq 2$, admits automorphisms with a dense orbit, 
and one can find such an automorphism contained in a subgroup $A$ of $\Aut (\T^n)\approx \GL(n,\Z)$ 
which is isomorphic to $\Z^d$ for $d\leq n-1$. This gives examples of   subgroups $A$ of $\Aut (\T^n) $ 
isomorphic to $\Z^d$ for $d\leq n-1$, acting with a dense orbit. Conversely every subgroup 
$A$ of $\Aut (\T^n) $ 
isomorphic to $\Z^d$ and acting with a dense orbit contains an (individual) automorphism which has a dense 
orbit; see \cite{Ber}; see also \cite{R10} for analogous results in a more general setting of 
automorphisms of general compact abelian groups. 

It is easy to see that we have a $\Z^2$-action on $\R$ with a dense orbit; the automorphisms defined 
respectively by
multiplication by $e^\alpha$ and $-e^\beta$, $\alpha , \beta>0$, generate such an action when 
$\alpha/\beta$ is 
irrational.  More generally, $G=\R^n\times \T^m$ admits a $\Z^d$-action with a dense orbit if any 
only if $m\neq 1$ and $d\geq (n+2)/2$; see \cite{D-Zd} for details.  Examples of nonabelian two-step 
nilpotent Lie groups $G$ admitting $\Z^2$-actions 
with a dense orbit are given in \cite{D-Zd}.

It may be mentioned here that Theorem~\ref{Zd} is extended in \cite{DSW} to actions on general locally 
compact groups $G$,
where it is concluded that under the analogous condition there exists a compact
normal subgroup $C$ such that the quotient $G/C$ is a (finite) product of locally compact fields
of characteristic zero. 

In line with the above it would be interesting to know about Lie groups admitting actions by nilpotent
or, more generally, solvable groups of automorphisms with a dense orbit.  

\subsection{Discrete groups with dense orbits}

Consider a connected Lie group $G$ with a Lie subgroup $H$ of  $\Aut (G)$ such  that the $H$-action on $G$ has an open dense orbit; e.g. 
$\R^n$, $n\geq 1$, and $H=\SL (n,\R)$  -- see \S \ref{denseorb} for more examples (recall also that such a $G$ is nilpotent). Let $g\in G$ be such that the $H$-orbit
 is open and dense in $G$. Let $L$ be the stabiliser of $g$ under the $H$-action, viz.
$L=\{\tau \in \Aut (G) \mid \tau (g)=g\}$. 
Then it is easy to see that for a subgroup $\Gamma$ of 
$H$ the $\Gamma$-action on $G$  has a  dense in $G$ if and only if the $L$-action on $H/\Gamma$
has a dense orbit; this phenomenon is known as ``duality" - see for instance \cite{BM}, \cite{D-encl}  for some details). 
In certain situations, such as when $\Gamma$ is a lattice in $H$ (viz. $H/\Gamma$ admits 
a finite measure invariant under the action of $H$ on the left), the question of 
whether the action of a subgroup has dense orbits on $H/\Gamma$ is amenable via techniques of ergodic theory.

For $G=\R^n$, $n\geq 2$, we have $\Aut (G)\approx \GL (n,\R)$, and there exist many discrete subgroups 
of the latter whose action on $\R^n$ admits dense orbits. Let $H=SL(n,\R)$ and 
$\Gamma$ be a lattice in $\SL (n,\R)$. We  choose $g$ as $e_1$ where $\{e_1,\dots ,e_n\}$ is the
coordinate basis of $\R^n$ and let $L$ be its stabiliser. It is known that this subgroup acts ergodically 
on $H/\Gamma$ and hence as noted above the $\Gamma$-action on $\R^n$ has a dense orbit. 
 This applies in particular to the subgroup $\SL(n,\Z)$,  consisting of integral unimodular matrices, which is indeed a lattice in $\SL (n,\R)$ (see \cite{MW}, Ch.\,7 or \cite{Rag}, Ch.\,10). We shall discuss more about the orbits of these
 in the next section. 

There are also natural examples of discrete subgroups of $\SL(n,\R)$ other than lattices which have dense orbits under the action on $\R^n$. For example if $n$ is even and $\Gamma$ is a lattice in the symplectic group $\Sp(n,\R)$ then 
it has dense orbits on $\R^n$; stronger statements analogous to those for lattices in $\SL (n,\R)$ are possible but 
we shall not go into the details. There are also other examples, arising from hyperbolic geometry. Let
$\Gamma$ be the fundamental group of a surface of constant negative curvature whose associated geodesic
flow is ergodic. Then $\Gamma$ may be viewed canonically as a subgroup of $\PSL (2,\R)=\SL (2,\R)/\{\pm I\}$,
where $I$ denotes the identity matrix, and if $\tilde \Gamma$ is the lift of $\Gamma$ in $\SL (2, \R)$, 
then the action of $\tilde \Gamma$ on $\R^2$ has dense orbits; the action is ergodic with respect to the 
Lebesgue measure  (see \cite{BM} for some details; see \S \ref{erg} for a discussion on ergodicity). Analogous examples can also be constructed 
in higher dimensions. 

In place of $\R^n$ one may consider other nilpotent connected Lie groups $G$ which admit an open dense 
orbit under the action of $\Aut (G)$. It would be interesting to know analogous results for discrete groups of automorphisms of  other connected nilpotent Lie groups, which however does not seem to be considered in the literature. 

\subsection{Orbit structure of actions of some discrete groups}

In most of the cases considered in the earlier subsections where we conclude existence of a dense orbit, not all orbits may be dense, and in general there is no good description possible of the ones that are not dense. However in certain cases a more complete description is possible. 

Expressing $\T^n$ as $\R^n/\Z^n$, $\Aut (\T^n)$ can be realised as $\GL(n,\Z)$, the group  
of all  $n\times n$ matrices with integer entries and determinant $\pm 1$. For $v\in \R^n$ whose 
coordinates with respect to the standard basis (generating $\Z^n$) are rational the orbits under
the $\Aut (\T^n)$-action are easily seen to be finite. It turns out, and not too hard to prove 
 (see \cite{DK}, for instance) that conversely for any $v$ at least one of whose coordinates is 
 irrational the $\Aut (\T^n)$-orbit is dense. 
 
 For the case of $\R^n$, $n\geq 1$ we noted above that if $\Gamma$ is a lattice in $\SL (n,\R)$ 
then there exist dense orbits, by duality and ergodicity considerations. In fact in this case it is possible 
to describe the dense orbits precisely. If $\SL (n,\R)/\Gamma$ is compact (viz. if $\Gamma$ is a 
``uniform" lattice) then the orbit of every non-zero point in $\R^n$ is dense in $\R^n$. When $\SL (n,\R)/\Gamma$ is noncompact (but 
a lattice) the set of points whose orbits are not dense is contained in a union of countably many lines 
(one-dimensional vector subspaces) in 
$\R^n$. In the case of $\SL(n,\Z)$ the exceptional  lines involved are precisely those passing 
through points in $\R^n$ with rational coordinates  (see \cite{D-Rag} for more general results in this direction). These results are consequences of the study of flows on homogeneous
spaces which has been a much studied topic in the recent decades, thanks to the work of Marina 
Ratner on invariant measures of unipotent flows. We shall not go into details on the topic
here; the interested reader is referred to the expository works \cite{D-encl} and \cite{KSS} and other references 
there, for exploration of the topic. The result in the special case of $\SL(n,\Z)$ recalled above was first proved in \cite{JSD} (see also \cite{DN} for a strengthening in another direction). 
 
 As in the case of the issues in the previous section, it would be interesting to know results analogous to the 
 above for discrete groups of automorphisms of other nilpotent connected Lie groups. 

\subsection{$\Aut (G)$-actions with few orbits}\label{forbits}

Observe that when $G=\R^n$, $n\in \N$, $\Aut (G)$ is $\GL(n,\R)$ and 
the action is transitive on the complement of the zero element.  Thus the action of $\Aut (G)$
on $G$ has only two orbits.  The question as to when there can be only finitely many, or countably many, orbits
has attracted attention, not only for Lie groups, but in the general context of 
locally compact groups; we shall indicate some of the results in that generality, giving 
references, but our focus shall be on Lie groups. 

\begin{theorem}\label{Strop}
Let $G$ be a  locally compact group. Suppose that the action of $\Aut (G)$ on $G$ has only countably many orbits.  Then the connected component $G^0$ of the identity in $G$ 
is a simply connected nilpotent Lie group. Moreover, the number of $\Aut (G)$-orbits in $G^0$ is finite and one of them is 
 an open orbit. 
\end{theorem}

The first statement in the theorem was proved in \cite{Str}. We note that for Lie groups it can be deduced from Theorem~\ref{denseorbit} (which is a generalised version of a
theorem from \cite{D-Sank}); since under the condition in the hypothesis at least one of the orbits 
has to be of positive measure Theorem~\ref{denseorbit} yields that the group is nilpotent, but on the other hand the 
condition also implies that there is no compact subgroup of positive dimension contained in the center,
so $G$ must in fact be a simply connected nilpotent Lie group. It is proved  in \cite{Str2} (Theorem~6.3) 
that under the condition in the hypothesis one of the orbits is open. Moreover, as $G$ is a simply connected
nilpotent Lie group, the $\Aut (G)$-orbits are in canonical one-one correspondence with $\Aut (\frak G)$-orbits on $\frak G$ via the exponential map, where $\frak G$ is the Lie algebra of $G$, and the latter being an action of a real algebraic group, the cardinality of the orbits can be countable only if  it is finite  and one of the 
orbits is open. 

Now let  $G$ be a connected  Lie group such that 
the $\Aut (G)$-action on $G$ has only finitely many orbits. 
Since the identity element is fixed, there  are
at least two orbits; the group is said to be {\it homogeneous}\index{homogeneous group} if there are only two orbits. 
The groups  $\R^n$, $n\geq 1$, are homogeneous, and they are also readily seen to be the only connected abelian Lie groups for which the $\Aut (G)$-action has only finitely many
orbits.  It turns out that $\R^n$, $n\geq 1$, are in fact the only connected 
 locally compact groups that are homogeneous (see \cite{Str98}, Theorem 6.4); all (not necessarily
connected) homogeneous locally compact groups have also been determined in \cite{Str98}; $K^n$, 
where $K$ is the field of $p$-adic numbers, with $p$ is a prime, or $\Q$ (with the discrete topology),
and $n\in \N$, are some of the other examples of homogeneous groups.

A locally compact group is said to be {\it almost homogeneous}\index{almost homogeneous group} if the $\Aut (G)$-action 
on $G$ has $3$ orbits.  The class of groups with the property has been studied in  \cite{Str} and \cite{Str2}. An almost homogeneous Lie group is a Heisenberg group\index{Heisenberg group}, namely a group 
defined on $V\oplus Z$, where $V$ and $Z$ are vector spaces, with the product 
defined by $(v,x)\cdot (w,y)=(v+w, x+y+\frac 12 \langle v,w\rangle)$, where $\langle \cdot, \cdot \rangle$
is an alternating non-degenerate bilinear form over $V$ with values in $Z$; (we note that here the group structure 
is viewed via identification with the corresponding Lie algebra, and the factor $\frac 12$ is introduced so that 
$\langle v,w\rangle$ is the Lie bracket of $v$ and $w$). Moreover, the pair of dimensions 
of the vector spaces $V$ and $Z$ are either of the form $(2n,1), (4n,2), (4n,3)$, with $n\in \N$, or
one of $(3,3), (6,6), (7,7), (8,5), (8,6)$ or $(8,7)$, and each of  these pairs determines an almost
homogeneous Heisenberg group.  Automorphism groups of these Heisenberg groups have been 
discussed in~\cite{Str99}.

It is easy to see that if $G$ is a $k$-step simply connected nilpotent Lie group 
then the number of $\Aut (G)$-orbits on $G$ is at least $k+1$. 
If $G$ is one of the Heisenberg groups as above then for any $r\geq 1$, the action of $\Aut (G^r)$ on $G^r$ 
has $2r+1$ orbits (cf. \cite{Str2}, Proposition~6.8). Thus there exist nilpotent Lie groups with arbitrarily large finite number of orbits under the action of the respective automorphism group. 
We refer the reader to \cite{Str2}, and other references 
there, for further details and also open problems on this theme, for general (not necessarily connected)
locally compact groups.

\section{Ergodic theory of actions of automorphism groups}\label{erg}

In this section we discuss various aspects of ergodic theory of 
actions on Lie groups by groups of automorphisms. 

\subsection{Preliminaries}

We begin by briefly recalling some definitions and conventions which will be followed throughout. 
By a measure we shall always mean a $\sigma$-finite measure. Given a measure space
 $(X, \frak M)$, an automorphism $\tau:X\to X$ is said measurable if $\tau^{-1}(E)\in \frak M$ 
 for all $E\in \frak M$; a measure $\mu$ on $(X,\frak M)$ is said to be invariant under a measurable 
 automorphism $\tau$ if $\mu (\tau^{-1}(E))=\mu (E)$ for all $E\in \frak M$, and it is said to be {\it quasi-invariant}\index{quasi-invariant measure} if, for $E\in \frak M$, $\mu (\tau^{-1}(E))=0$ if and only if $\mu (E)=0$; a measure is said to be
 invariant or, respectively, quasi-invariant, under a group of 
measurable  automorphisms if it has the property with respect to the action of each of the automorphisms
from the group. A measure which is invariant or quasi-invariant with respect to an action is 
said to be {\it ergodic}\index{ergodic action} with respect to the action if there are no two disjoint measurable
subsets  invariant under the action, each with positive measure. 

Two measures on a measurable 
space $X$ are said to be equivalent if they have the same sets of measure $0$.  It is easy to see 
that given a measure which is quasi-invariant under an action  there exists a finite measure equivalent to it,
which is also quasi-invariant; in particular, given an infinite invariant measure there exists a finite quasi-invariant measure equivalent to it (which may however not be invariant under the action).  

A measure $\mu$ on $G$ is called  a probability measure if $\mu (G)=1$. 
We denote by $P(G)$ the space of probability measures on $G$. 
For actions of a locally compact second countable group every quasi-invariant probability measure can be 
``decomposed as a 
continuous sum" of ergodic quasi-invariant probability measures in a canonical way, and hence it suffices in many respects to understand 
the ergodic quasi-invariant measures (see, for instance, \cite{MW}, Theorem~14.4.3).

In our context the measure space structure will always be with respect to the Borel $\sigma$-algebra
of the topological space in question.

\subsection{Finite invariant measures}

When $G$ is the torus $\T^n$, $n\geq 2$, the Haar measure is invariant under the action of $\Aut (G)$, 
which is an infinite discrete group. More generally, if $G$ is a connected Lie group such that the center 
contains a torus of positive dimension, then the Haar measure of the maximal torus of the center, viewed 
canonically as a measure on $G$, is invariant under $\Aut (G)$, which can be a group with infinitely 
many connected components (see \S \ref{aa}). It turns out however that when we restrict to almost 
algebraic subgroups of $\Aut (G)$ the situation is quite in contrast, as will be seen in Theorem~\ref{finiteinv} below.

Let $G$ be a connected Lie group. 
For $\mu\in P(G) $ we denote by $\supp \mu$ the support of $\mu$, namely the smallest closed subset 
of $G$ whose complement is of $\mu$-measure $0$.  Let   
$$I(\mu)=\{\tau \mid \mu \hbox{ \rm is invariant under the action of } \tau\},$$  
and 
$$J(\mu)=\{\tau \in \Aut (G) \mid \tau (g)=g \hbox{ \rm for all } g\in \supp \mu\}.$$
Then it can be seen that $I(\mu)$ and $J(\mu)$ are both closed subgroups of $\Aut (G)$ and   $J(\mu)$ 
is a normal subgroup of $I(\mu)$. It turns out that when $G$ has no compact subgroup of positive 
dimension contained in the center, the quotient $I(\mu)/J(\mu)$ is compact. In fact we have
the following. 

\begin{theorem}\label{finiteinv}{\rm (cf. \cite{D-inv})}
Let $G$ be a connected Lie group. 
Let $\frak A$ be an almost algebraic subgroup of $\Aut (G)$. Then for any $\mu \in P(G)$, 
$(I(\mu) \cap \frak A)/ (J(\mu) \cap \frak A)$ is compact. 

\end{theorem}

Consider a connected Lie group  $G$  such that $\Aut (G)$ 
is almost algebraic (see Theorem~\ref{denseorbit}) then $I(\mu)/J(\mu)$ is compact and in turn 
for any $x\in \supp \mu$ the $I(\mu)$-orbit of $x$  in $G$ is compact; thus $\supp \mu$ can 
be expressed as a disjoint union of compact orbits of $I(\mu)$.

We recall that if $G$ is a compact Lie group then the Haar measure of $G$ is a finite measure 
invariant under all automorphisms. On the other hand when $G=\R^n$ for some $n\geq 1$, the only finite 
invariant measure is the point mass at~$0$ (see Corollary~\ref{cor:4.2} below). The general situation is, in a sense, a mix of the two
kinds of situations.

From Theorem~\ref{finiteinv}  one can  deduce the following, by going modulo the maximal 
compact subgroup contained in the center and applying Theorem~\ref{denseorbit} on the quotient. 

\begin{corollary}\label{cor:4.2}
Let $G$ be a connected Lie group and  $H$ be a subgroup of $\Aut (G)$.  Let  $\mu$ be a 
finite $H$-invariant measure on $G$. Then for any $g\in \supp \mu$ the $H$-orbit of $g$ is contained
in a compact subset of $G$. In particular if $H$ does not have an orbit other than that of the
identity which is  contained in a compact subset, then the point mass at the identity is the only 
$H$-invariant probability measure on $G$. 
\end{corollary} 

For $G=\R^n$, $n\geq 1$, Therorem~\ref{finiteinv} implies in particular the following result, 
which was proved earlier in  \cite{D-quasi}. 

\begin{theorem}
For $\mu \in P(\R^n)$, $I(\mu)$ 
is an algebraic subgroup of $\Aut (\R^n)= \GL(n,\R)$.

\end{theorem}

Analogous results are also proved in \cite{D-quasi} for projective transformations. The approach 
involves the study of non-wandering points of the transformations, which has been discussed in
further detail in \cite{D-04}. 

\subsection{Convolution powers of invariant measures}

We recall here the following result concerning convolution powers of probability measures invariant under the action of a compact subgroup of positive dimension. For a probability measure $\mu$ and $k\in \N$ we denote 
by $\mu^k$ the $k$-fold convolution power of $\mu$.

\begin{theorem}\label{K-inv}
Let  $\mu$ be a probability measure on $\R^n$, $n\geq 1$. Suppose that $\mu$
is  invariant under the action of a compact 
connected subgroup $K$ of $\Aut  (G)=\GL (n,\R)$ of positive dimension. Then one of the following holds:

i) $\mu^n$ is absolutely continuous with respect to the Lebesgue measure on $\R^n$;

ii) there exists an affine subspace $W$ of $\R^n$ such that $\mu (W)>0$. 

\end{theorem}

This is a variation of  Theorem~3.2 in \cite{DGS} whose proof can be read off from that of the latter; it 
can be seen  that Condition (ii) above fails to hold if the  $K$-action 
has no nonzero fixed point and there is no proper vector subspace $U$ with $\mu (U)>0$ (the subspace $U$ 
can  also be stipulated to be invariant), so under these assumptions $\mu^n$ is absolutely continuous 
with  respect to the Lebesgue measure; this is the formulation of Theorem~3.2 in \cite{DGS}. 

It would be interesting to have a suitable analogue of the above theorem for actions on a general connected Lie 
group $G$ by compact subgroups of positive dimension in $\Aut (G)$, (with the Lebesgue 
measure replaced by the Haar measure of $G$). 

\subsection{Infinite invariant measures}

The Lebesgue measure on $\R^n$ is invariant under the 
action of a large subgroup of $\Aut (G) =\GL(n,\R)$, namely $\SL(n,\R)$, the special linear group, but not 
the whole of $\Aut (G)$.  On the other hand if $G$ is a connected semisimple Lie group then $\Inn (G)$ 
is of finite index in $\Aut (G)$, and as $G$ is unimodular it follows that the Haar measure is invariant under 
the action of $\Aut (G)$. 
In many Lie groups of common occurrence  it is readily possible to determine whether the Haar measure is invariant under all automorphisms, but there does not seem to be
any convenient  criterion in the literature to test it. 

Let $G$ be a connected Lie group and let $S$ denote the subgroup which is the center of $\Aff^0(G)$, 
when $G$ is viewed canonically as a subgroup of $\Aff (G)$ (see~\S~\ref{linrep}). Recall that by 
Theorem~\ref{locallyclosed} all orbits of $\Aut^0(G)$, and more generally of any almost algebraic 
subgroup $A$ of $\Aut (G)$, on $G/S$ are locally closed. By a theorem of Effros~\cite{Ef} 
this implies that for every measure on $G$ which is quasi-invariant and ergodic under the action of an 
almost algebraic subgroup $A$ of $\Aut (G)$, there exists an $A$-orbit $\frak O$ such that the complement
of $\frak O$ has measure $0$ (that is, the measure is ``supported" on  $\frak O$, except that the latter
is a locally closed subset which may not be closed);   in particular this applies to any infinite measure 
invariant under the action of an almost algebraic subgroup $A$ of $\Aut (G)$ as above. 

\subsection{Quasi-invariant measure and ergodicity}

  We recall that for an 
action of  a group on a locally compact
second countable space, given an ergodic quasi-invariant measure $\mu$, for almost all points in $\supp \mu$,
the orbit is dense in $\supp \mu$. Together with the results on existence of dense orbits 
in~\S\,3 this yields the following:

\begin{theorem}
Let $G$ be a connected Lie group and let $\lambda$ be a Haar measure of $G$. Then $\lambda$ 
is quasi-invariant   under the action of $\Aut (G)$. Moreover  the
following holds:

i) if $\lambda$  is ergodic with respect to the  
$\Aut (G)$-action on $G$ then $G$ is a nilpotent Lie group. 

ii) if $\lambda$ is ergodic  under the action of an abelian subgroup of $\Aut (G)$ then $G$ is a two-step 
nilpotent Lie group, with $\overline {[G,G]}$ compact. 

iii) if $\lambda$ is ergodic under the action of some $\alpha \in \Aut (G)$ then $G$ is a torus. 

iv)  if $\lambda$ is ergodic under the action of an almost algebraic subgroup $A$ of $\Aut (G)$ then the 
$A$-action on $G/C$, where $C$ is the maximal compact subgroup of $G$, has an open dense orbit. 

\end{theorem}

The first three statements are straightforward consequences of results from~\S\,3 on the existence of a dense orbit and the observation above.  The last 
assertion follows from the fact that $\{\bar \alpha \mid \alpha \in A\} $ is an almost algebraic subgroup
of $\Aut (\frak G')$, where $\frak G'$ is the Lie algebra of $G/C$, and hence its orbits on $\frak G'$ 
are locally closed; as $G/C$ is a simply connected nilpotent Lie group the exponential map is a 
homeomorphism and hence  the preceding conclusion implies that the $A$-orbits on $G/C$ are locally closed, and 
in particular an orbit which is dense is also open in $G/C$.

\subsection{Stabilisers of actions of Lie groups}\label{4.4}

Let $G$ be a connected Lie group and consider a (Borel-measurable) action of $G$ on a standard Borel space $X$. By the stabiliser of a point $x\in X$ we mean the subgroup $\{g\in G\mid gx=x\}$, and we denote it by
$G_x$. Each  $G_x$ is a closed subgroup of $G$ (see \cite{V-geom}, Corollary~8.8). 
When the action is 
transitive, namely when the whole of $X$ is a single orbit, then the stabilisers of any two points in $X$ are conjugate to each other. One may wonder to what extent this generalises to a general ergodic action, with respect
to a measure which is invariant or quasi-invariant under the action of $G$. In this respect the following 
is known from \cite{D-conj}; the results in \cite{D-conj} strengthen those from \cite{Gol} in the case of Lie 
groups, while in \cite{Gol} the issue of conjugacy of the stabilisers is considered in the wider framework of all locally compact groups. 

\begin{theorem}\label{clsub} 
Let $G$ be a connected Lie group acting on a standard Borel space $X$. Let $\mu$ be a
measure on $X$ which is quasi-invariant and ergodic with respect to the action of $G$. Suppose also that 
for all $x\in X$ the stabiliser $G_x$ of $x$ has only finitely many connected components. Then there 
exists  a subset $N$ of $X$ with $\mu(N)=0$ such that for $x,y\in X\backslash N$ there exists $\alpha\in 
\Aut^0(G)$ such that $\alpha (G_x)=G_y$; in particular $G_x$ and $G_y$ are topologically isomorphic to 
each other. 
\end{theorem}

This is a variation of Corollary 5.2 in \cite{D-conj}, where the ergodicity condition involved is formulated 
in terms of 
a $\sigma$-ideal of null sets, rather than a quasi-invariant measure. Without the assumption of $G_x$ 
having only finitely many connected components it is proved that for $x,y$ in the complement of a null 
set $N$, the connected components of the identity in $G_x$ and $G_y$ are topologically isomorphic to each 
other (\cite{D-conj}, Corollary~5.7). The proofs of these results are based on consideration of the action of 
$\Aut (G)$ on the space of all closed subgroups of $G$, the latter being equipped with the Fell topology, and the orbits 
or stabilisers $G_x$, $x\in X$, under the action. The argument depends on the fact $\Aut^0(G)$ is 
an almost algebraic subgroup of $\GL(\frak G)$ containing $\Inn (G)$, $\frak G$ being the Lie algebra of $G$ (see Theorem~\ref{thm:1.2}(ii)); the proof shows that automorphisms $\alpha$ as in the conclusion of the 
theorem may in fact be chosen to be from the smallest almost algebraic subgroup containing $\Inn (G)$, and in particular if $\Inn (G)$ is an almost algebraic subgroup of $\GL(\frak G)$ then stabilisers of almost all points are conjugates to each other.  Other conditions under which the stabilisers of almost all points may be concluded to be actually conjugate in $G$, and also examples for which it does not hold,  are discussed in \cite{D-conj}. 
We shall not go into the details of these here. 

\section{Some aspects of  dynamics of $\Aut (G)$-actions}

In this section we discuss certain results on dynamics in which some features of the actions of $\Aut (G)$ on $G$ described in earlier sections play a role. 

\subsection{Stabilisers of continuous actions}\label{5.1}

Let $G$ be a connected Lie group and consider a continuous action of $G$ on a compact Hausdorff space $X$. Here we recall a topological analogue of Theorem~\ref{clsub}. As before we denote by $G_x$ 
the stabiliser of $x$ under the action in question and we denote by $G_x^0$ the connected component 
of the identity in $G_x$. We recall that an action is said to be {\it minimal}\index{minimal action} if there 
is no proper nonempty closed subset invariant under the action.  The following is proved in 
\cite{D-stab}, Proposition~3.1.

\begin{theorem}
Let $G$ be a connected Lie group acting continuously on a locally compact Hausdorff space $X$. Suppose that 
there exists $x\in X$ such that the orbit of $x$ is dense in $X$. Then there exists an open dense subset $Y$ of $X$ such that for all $y\in Y$ there exists $\alpha \in \Aut^0(G)$ such that $\alpha (G_x^0)=G_y^0$. If the action is minimal then $G_x^0$, $x\in X$, are Lie isomorphic to each other. 
\end{theorem}

Further conditions which ensure $G_x^0$, $x\in X$, being conjugate in $G$ are discussed in \cite{D-stab}. 
The results of \cite{D-stab} generalise earlier results of C.\,C. Moore and G. Stuck proved for special classes of Lie groups (see \cite{D-stab} for details). 

\subsection{Anosov automorphisms}

Anosov diffeomorphisms have played an important role in the study of differentiable dynamical systems; we shall not go 
into much detail here on the issue -- the reader is referred to Smale's  expository article \cite{Sm} for a  perspective on the topic. Hyperbolic automorphisms of tori $\T^n$, $n\geq 2$, 
namely automorphisms $\alpha$ such that $d\alpha $ has no eigenvalue (including complex) of absolute value $1$, serve as the simplest examples of Anosov automorphisms. In \cite{Sm} Smale described an example, due
to A. Borel, of a non-toral compact nilmanifold admitting Anosov automorphisms; a nilmanifold is a homogeneous space of the form $N/\Gamma$ where $N$ is a simply connected nilpotent Lie group and $\Gamma$ is discrete subgroup, and an Anosov automorphism of $N/\Gamma$ is the quotient on $N/\Gamma$ of an automorphism $\alpha$ of $N$ such that $\alpha (\Gamma)=\Gamma$ and $d\alpha$ has no eigenvalue of absolute value $1$. Such a system can have a nontrivial finite group of symmetries and factoring through 
them leads to some further examples on what are called infra-nilmanifolds, known as Anosov automorphisms of infra-nilmanifolds. Smale conjectured in \cite{Sm} that all Anosov 
diffeomorphisms are topologically equivalent to Anosov automorphisms on infra-nilmanifolds. 

A broader class of Anosov automorphisms of nilmanifolds, which includes also the example of Borel referred to above, was introduced by L.~Auslander and J.~Scheuneman \cite{AS}. Their approach involves analysing $\Aut (N)$  
for certain simply connected nilpotent Lie groups $N$ to produce examples of Anosov automorphisms of 
compact homogeneous spaces of $N/\Gamma$ for certain discrete subgroups $\Gamma$. In \cite{D-nil} the 
approach in \cite{AS} was extended using some results from the theory of algebraic groups and arithmetic 
subgroups, and some new examples of Anosov automorphisms were constructed.  A large class of examples of Anosov automorphisms were constructed in \cite{DM-a} by associating  nilpotent Lie groups to graphs and studying their automorphism group in relation to the graph. Study of the automorphism groups of nilpotent Lie groups has also been applied to construct examples of nilmanifolds that can not admit Anosov automorphisms (see \cite{D-nil} and \cite{D-step}). 

Subsequently Anosov automorphisms of nilmanifolds have been constructed via other approaches; the
reader is referred to \cite{G} and \cite{DV} and various references cited there for further details. There 
has however been no characterisation of nilmanifolds admitting Anosov automorphisms, and it may be hoped that further study of the automorphism groups of nilpotent Lie groups may throw more light on the issue.

\subsection{Distal actions}

An action of a group $H$ on a topological space $X$ is said to be {\it distal}\index{distal action} if for any pair of distinct points 
$x,y$ in $X$ the closure of the $H$-orbit $\{(gx,gy)\mid g\in G\}$ of $(x,y)$ under the componentwise action  of $G$ on the 
cartesian product space $X\times X$
does not contain a point on the ``diagonal'', namely of the form $(z,z)$ for any $z\in X$. The action on a locally compact group $G$ 
by a group $H$ of  automorphisms of $G$ is distal if and only if under the $H$-action on $G$ the closure the orbit of any nontrivial element $g$ in $G$ does not contain the identity element. 

Distality of actions is a classical topic, initiated by D. Hilbert, but the early studies were limited to actions on compact spaces. 
The question of distality of actions on $\R^n$, $n\geq 1$, by groups of linear transformations was initiated by 
C.\,C. Moore \cite{M} and was strengthened by Conze and Guivarc'h \cite{CG}; see also  H.~Abels~\cite{A}, 
where the results are extended to actions by affine transformations. It is proved that the action of a group $H$ of linear automorphisms of $\R^n$ is distal if and only if the action of each $h\in H$ (viz. of the cyclic 
subgroup generated by it) is distal, and that it holds if and only if all (possibly complex) eigenvalues of $h$ are 
of absolute value $1$ (see \cite{CG}, \cite{A}). 

It is proved in \cite{A-group} that the action of a group $H$ of automorphisms of a connected Lie group $G$ is distal if and only if the 
associated action of $H$ on the Lie algebra $\frak G$ of $G$ is distal (to which the above characterisations 
would apply). 

For a subgroup $H$ of $\GL (n,\R)$, the $H$-action  
on $\R^n$  is distal if and only if the action of its algebraic hull (Zariski closure) in $\GL (n,\R)$ is distal, and 
the action of an algebraic subgroup $H$ of $\GL (n,\R)$ on $\R^n$ is distal if and only if the unipotent 
elements in $H$ form a closed subgroup $U$ (which would necessarily be the unipotent radical) such that $H/U$  is compact. (see \cite{A}, Corollaries~2.3 and~2.5). Via the above correspondence these results can be applied also to actions of 
groups of automorphisms of a general connected Lie group $G$.

An action of a group $H$ on a space $X$ is called  {\it MOC} (short for ``minimal orbit closure")\index{MOC action} if 
the closures of all orbits are minimal sets (viz. contain no 
proper nonempty invariant closed subsets). When $X$ is compact, the MOC condition is equivalent to distality. 
For the action on a topological group $G$ by a group of automorphisms MOC implies distality.  
The converse is known in various cases, including 
for actions on connected Lie groups $G$ 
by groups of automorphisms  (see \cite{A-group}). The reader is also referred to \cite{A-genr}, \cite{Sh-12} 
and \cite{RS2} for some generalisations of this as well as some of the other properties discussed above to more general 
 locally compact groups).

A connected Lie group $G$ is said to be of {\it type $\cal R$}\index{Lie group of type $\cal R$} if  for all $g\in G$ all (possibly complex) eigenvalues of $\Ad \,g$ are of absolute value $1$. Thus in the light of the results noted above $G$ is of type $\cal R$ if and only if the action of $\Inn (G)$ on $G$  is distal. This condition is also equivalent to  $G$ having polynomial growth\index{group of polynomial growth}, viz. for any compact neighbourhood $V$ of the identity the 
Haar measure of $V^n$ grows at most polynomially in $n$; see \cite{Jen} and also \cite{Gui} and \cite{Los} for related general results. 

\subsection{Expansive actions} 

A homeomorphism $\varphi$ of a compact metric space $(X,d)$ is said to be {\it expansive}\index{expansive action} if there exists $\epsilon >0$ such that  for any  pair of distinct points $x,y$ there exists an integer $n$ such that $d(\varphi^n(x), \varphi^n(y))>\epsilon$; the notion may also be defined in terms of a uniformity in place of a metric. An automorphism $\tau$ of a topological group  $G$ is expansive if there exists a neighbourhood $V$ of the identity such that for any nontrivial element $g$ in $G$ there exists an integer $n$ such that $\tau^n(g)\notin V$. More generally the action of a group $\Gamma$ of automorphisms of a topological group $G$ is said to be expansive if there exists a neighbourhood $V$ of the identity in $G$ such that 
for any nontrivial $g$ in $G$ there exists a $\gamma \in \Gamma$ such that $\gamma (g) \notin V$.  

A compact connected topological group admits a group of automorphisms acting expansively only if it is abelian (cf. \cite{Lam}) and finite dimensional (cf. \cite{Law}). For compact abelian groups the expansiveness condition on actions of automorphism groups has been extensively studied using techniques of commutative algebra (see \cite{Sch} for details). 

If $G$ is a connected Lie group with Lie algebra $\frak G$, the action of a subgroup $\Gamma$ of $\Aut (G)$ 
on $G$ is expansive if and only if for the induced action of $\Gamma$ on $\frak G$, for any nonzero $\xi \in \frak G$ the $\Gamma$-orbit of $\xi$ is unbounded in $\frak G$  (cf. \cite{Bh}, where the issue is also considered for semigroups of endomorphisms). It is also deduced in \cite{Bh} that if $\Gamma$ is a virtually nilpotent Lie group (viz. with a nilpotent subgroup of finite index) of $\Aut (G)$ whose action on $G$ is expansive, then $\Gamma$ 
contains an element acting expansively on $G$.  

\section{Orbits on the space of probability measures}

The action of $\Aut (G)$ on $G$, where $G$ is a Lie group, induces an action of $\Aut (G)$ on the space
of probability measures on $G$ (see below for details). This action plays an important role in many contexts. This section is devoted 
to recalling various results about the action and their applications. 

\subsection{Preliminaries}

Let $G$ be a connected Lie group and as before let $P(G)$ denote the space of probability measures on 
$G$. We consider $P(G)$ equipped with the weak$^*$ topology with respect to the space of bounded continuous functions; we note that the topology is metrizable,  
and a sequence $\{\mu_j\}$ in $P(G)$ converges to $\mu\in P(G)$ if and only if $\int_Gf d\mu_j\to \int_G
fd\mu$, as $j\to \infty$, for all bounded continuous functions $f$ on $G$. We recall that the action of $\Aut (G)$ on $G$ induces an 
action on $P(G)$, defined, for  $\tau \in \Aut (G)$ and $\mu \in P(G)$, by $\tau (\mu)(E)=\mu (\tau^{-1}(E))$ 
for all Borel subsets $E$ of $G$. 

For $\lambda , \mu\in P(G)$, $\lambda * \mu$ denotes the convolution product of $\lambda$ and $\mu$, and for any $\mu \in P(G)$ and $k\in \N$ we denote by $\mu^k$ the $k$-fold convolution product $\mu * \cdots * \mu$. 

Behaviour of orbits of probability measures under actions of various subgroups is of considerable interest in various contexts. One of the issues is to understand conditions under which  orbits of automorphism groups 
$A\subset \Aut (G)$ are locally 
closed, namely open in their closures. We recall here that by a result of Effros \cite{Ef} this condition 
is equivalent to a variety of other ``smoothness" conditions for the action, including that the orbits map being
an open (quotient) map onto  its image (the latter being considered with respect to the induced topology from $P(G)$). 

As before for a connected Lie group $G$ we view $\Aut (G)$ as a subgroup of $\GL (\frak G)$, where $\frak G$ is the Lie algebra of $G$, and a subgroup $A$ of $\Aut (G)$ is said to be almost algebraic if it is an almost algebraic subgroup of $\GL(\frak G)$ (see \S~\ref{algebraic}). 

\begin{theorem}\label{measureorbits}({\rm cf. \cite{DR}, Theorem~3.3})
Let $G$ be a connected Lie group and $A$ be an almost algebraic subgroup of $\Aut (G)$. Let $C$ be 
the maximal compact subgroup contained in the center of $G$. Suppose that for any $g\in G$, 
$\{g^{-1}\tau (g)\mid \tau \in A\}\cap C$  is finite. Then for any $\mu \in P(G)$ the $A$-orbit $\{\tau (\mu)\mid 
\tau \in A\}$ is open in its closure in $P(G)$. Moreover, if $A$ consists of  unipotent elements in 
$\GL (\frak G)$ then the $A$-orbit is closed in $P(G)$. 

\end{theorem}

The theorem implies in particular that when $G$ has no compact subgroup of positive dimension 
contained in its center, for the action of any almost algebraic subgroup $A$ of $\Aut (G)$ (which includes 
also the whole of $\Aut (G)$ in this case) the orbits of $A$ are locally closed. 

\subsection{Factor compactness and concentration functions}

Let $\mu\in P(G)$. We denote by $G(\mu)$ the smallest closed subgroup of $G$ containing $\supp \mu$, 
the support of $\mu$, by $N(\mu)$ the normaliser of $G(\mu)$ in $G$ and by $Z(\mu)$ the centraliser of 
$\supp \mu$, namely $\{g\in G\mid gx=xg \hbox{ \rm for all }x \in \supp \mu\}$. 

A $\lambda \in P(G)$ is called a {\it factor}  of $\mu$ if there exists $\nu\in P(G)$ 
such that $\mu =\lambda *\nu=\nu *\lambda$. It is known that any factor of $\mu$ has its support 
contained in $N(\mu)$,  more specifically in a coset of $G(\mu)$ contained in $N(\mu)$ (see \cite{DM-f}, 
Proposition~1.1 for the first assertion; the second is an easy consequence of the first). 

A question of interest, which is not yet  understood 
in full generality,
is whether given $\mu \in P(G)$ and a sequence $\{\lambda_j\}$ of its factors there exists a 
sequence $\{z_j\}$ in $Z(\mu)$ such that $\{\lambda_jz_j\}$ is relatively compact 
in $P(G)$ (in place of $Z(\mu)$ one may also allow in this respect the subgroup $\{g\in G\mid
g\mu g^{-1}=\mu\}$, which can be bigger, but the distinction turns out to be a rather technical issue, and we shall
not concern ourselves with it here). It is known that given a sequence of factors $\{\lambda_j\}$ there exists a sequence 
$\{x_j\}$ in $N(\mu)$ such that $\{x_j\mu x_j^{-1}\}$, $\{x_j^{-1}\mu x_j\}$ and $\{x_j\lambda_j\}$ are relatively compact \cite{DM-f}, \cite{DM}, and it suffices to show that for such a sequence $\{x_j\}$ the sequence 
of cosets $\{x_jZ(\mu)\}$ is relatively 
compact in $G/Z(\mu)$.  It is known that this is true when  $G$ is an almost algebraic Lie group (and also under 
some weaker, somewhat technical, conditions) (see \cite{DM-f}, \cite{DR}); the proofs however fall back on reducing
the question to vector space situation and a more direct approach connecting the question to 
Theorem~\ref{measureorbits} would be desirable. 

A similar issue, or rather a more general one, arises in the study of the decay of concentration 
functions\index{concentration function} of $\mu^n$ as $n\to \infty$, where one would like to know under what conditions a sequence
of the form $\{x_j\mu x_j^{-1}\}$ is relatively compact. We shall however not go into the details of the 
concept or the 
results about it here; a treatment of the topic in this perspective may be found in  \cite{D-Par} and \cite{DS-conc}; for a complete resolution of the problem itself, which however involves
a somewhat different approach, the reader is referred to  \cite{Jaw} and other references cited there.  

\subsection{Tortrat groups}
A locally compact group is called a {\it Tortrat group} if for any $\mu \in P(G)$ which is not idempotent
(viz. such that $\mu^2\neq \mu$) the closure of $\{g\mu g^{-1} \mid g\in G\}$ in $P(G)$ does not contain any 
idempotent measure. We note that the issue concerns the closure of the orbits of $\mu$ in $P(G)$ 
under the action of the subgroup $\Inn (G)$; the latter being contained in $\Aut^0(G)$ which is an 
almost algebraic subgroup is useful in this respect. 

Recall (see \S\,5.3) that a connected Lie group is said to be type {\it $\cal R$} if for all 
$g\in G$ all eigenvalues of $\Ad g$ are of absolute value $1$. It was shown in 
\cite{DR2} that for a Lie group the two properties are equivalent:

\begin{theorem}\label{tortrat}{\rm (cf. \cite{DR2})}
A connected Lie group is a Tortrat group if and only if it is of {\it type $\cal R$}. 

\end{theorem}

The result implies also that a connected locally compact group is a Tortrat group if and only if it 
is of polynomial growth (see \cite{DR2} for details and the related references). 

\subsection{Convergence of types}

Let $\{\lambda_j\}$ and $\{\mu_j\}$ be two sequences in $P(G)$ converging to $\lambda$ and $\mu$ 
respectively. Suppose further that there exists a sequence $\{\tau_j\}$ in $\Aut (G)$ such that $\tau_j(\lambda_j)=\mu_j$ for all $j$. A question, arising in various contexts in the theory of probability 
measures on groups, is under what (further) conditions on $\lambda$ 
and $\mu$ can we conclude that there exists a $\tau \in \Aut (G)$ such that $\tau (\lambda)=\mu$; 
the reader can readily convince herself/himself that further conditions are indeed called for. When there
exists a $\tau$ as above we say that {\it convergence of types} holds. Though the term ``type of a measure" does not seem to make an appearance in literature freely, 
implicit in the terminology above is the idea that $\{\tau (\mu)\mid \tau \in \Aut (G)\}$ constitutes the 
``type of $\mu$", or that $\mu$ and $\tau (\mu)$ are of the same type for any $\tau \in \Aut (G)$, and the 
question is if you  have pair of sequences of measures with the corresponding measures of the same type, converging to a pair of measures, 
under what conditions can we conclude the limits to be of the same type. 

We now recall some results in this respect; for reasons of simplicity of exposition we shall not strive for full generality (see \cite{D-inv} for more details). As in \S\,2 let $\Aff (G)$ 
be the group of affine automorphisms, and let $\rho :\Aff (G) \to \GL(V)$ be the linearising representation (see 
\S \ref{linrep});
recall that $\rho$ is defined over $V=\wedge^a \frak B$, where $\frak B$ is the Lie algebra of $\Aff (G)$ and
$a$ is the dimension of $\Aut (G)$. 

We say that $\mu \in P(G)$ is {\it $\rho$-full} if there does not exist any proper subspace $U$ of $V$ 
that is invariant under $\rho (g)$ for all $g\in \supp \mu$; though this condition is rather technical, as it 
involves the linearising representation, it is shown in \cite{D-inv} that for various classes of groups it holds 
under simpler conditions; for instance if $G$ is an almost algebraic subgroup of 
$\GL(n, \R)$ then the condition holds for any $\mu$ for which $\supp \mu$ is not contained in 
a  proper almost algebraic subgroup of $G$. The following is a result which is midway between 
the rather technical Theorem~1.5 and the  specialised result in Theorem~1.6 from \cite{D-inv}, whose 
proof can be read off from the proof of Theorem~1.6. 

\begin{theorem}\label{convoftypes}
Let $G$ be a connected Lie group. Let $\{\lambda_j\}$ and $\{\mu_j\}$ be sequences in $P(G)$ converging to $\lambda$ and $\mu$ 
respectively. Suppose that there exists a sequence $\{\tau_j\}$ in $\Aut^0 (G)$ such that $\mu_j=\tau_j(\lambda_j)$ for all $j$ and that $\lambda$ and $\mu$ are $\rho$-full. Then there exist sequences
$\{\theta_j\}$ and $\{\sigma_j\}$ in $\Aut (G)$ such that $\{\theta_j\}$ is contained in a compact subset 
of $\Aut (G)$, $\{\sigma_j\}$ are isotropic shear automorphisms, and 
$\tau_j=\theta_j\sigma_j$ for all $j$. If moreover $\supp \lambda$ is not contained in any proper closed 
normal subgroup $M$ of $G$ such that $G/M$ is a vector space of positive dimension, then $\{\sigma_j\}$ is also relatively 
compact. 

\end{theorem}

We note that when $\{\tau_j\}$ as in the theorem is concluded to be relatively compact, convergence of 
types holds (for the given sequences); with the notation as in the theorem, if $\tau$ is an accumulation point of $\{\tau_j\}$ then $\tau (\lambda) =\mu$. In the light of the conclusion of the theorem it remains mainly to understand the asymptotic behaviour under sequences of shear automorphisms.  

A special case of interest is when $\{\lambda_j\}$ are all equal, say $\lambda$. For this case we recall also
the following result, for shear automorphisms, proved in \cite{DGS} (Theorem~4.3 there);  the result played an important role in the proof of the main theorem  there concerning embeddability in a continuous one-parameter semigroup, for a class of infinitely divisible probability measures; we shall not go into the details of these concepts here. It may be noted that one of the conditions in the hypothesis is as in statement~(i) of Theorem~\ref{K-inv}. 

\begin{theorem}
Let $G$ be a connected Lie group and $T$ be a torus contained in the center of $G$. Let $H$ be a closed
normal subgroup of $G$ such that $G/H$ is topologically isomorphic to $\R^n$ for some $n\geq 1$. Let 
$W=G/H$ and  $\theta :G\to W$ be the canonical quotient homomorphism. Let $\lambda \in P(G)$ be 
such that $\theta (\mu)^n$ 
is absolutely continuous with respect to the Lebesgue measure on $W$. Let $\{\varphi_j\}$ be a 
sequence of homomorphisms of $W$ into $T$ and for all $j$ let $\tau_j\in \Aut (G)$ be the shear automorphism 
of $G$ corresponding to $\varphi_j$. Suppose that the sequence $\{\tau_j(\lambda)\}$ converges to a measure 
of the form $\tau (\lambda)$ for some $\tau \in \Aut (G)$ fixing $T$ pointwise. Then there exists a sequence 
$\{\sigma_j\}$ in $\Aut (G)$ such that $\sigma_j(\mu)=\mu$ and $\{\tau_j\sigma_j\}$ is relatively compact. 
\end{theorem}

\vskip1cm

\begin{flushleft}
S.G. Dani\\
Department of Mathematics\\
Indian Institute of Technology Bombay\\
Powai, Mumbai 400005\\
India

\vskip3pt
E-mail: {\tt sdani@math.iitb.ac.in}

\end{flushleft}

\printindex

\end{document}